\documentclass{elsart3-1-r}
\usepackage{amssymb}
\usepackage[english,francais]{babel}

\newtheorem{theorem}{Theorem}
\newtheorem{lemma}[theorem]{Lemma}
\newtheorem{e-proposition}[theorem]{Proposition}
\newtheorem{corollary}[theorem]{Corollary}
\newtheorem{e-definition}[theorem]{Definition\rm}
\newtheorem{remark}{\it Remark\/}

\newtheorem{theoreme}{Th\'eor\`eme}
\newtheorem{lemme}[theoreme]{Lemme}

\newtheorem{corollaire}[theoreme]{Corollaire}

\renewcommand{\theequation}{\arabic{equation}}
\def\og{\leavevmode\raise.3ex\hbox{$\scriptscriptstyle\langle\!\langle$~}}
\def\fg{\leavevmode\raise.3ex\hbox{~$\!\scriptscriptstyle\,\rangle\!\rangle$}}

\newcommand{\iomu}[1]{\int #1\,d\mu}
\newcommand{\iomuz}[1]{\int #1\,d\nu}
\newcommand{\proof}{\noindent{\sl Proof.\/}\ }
\newcommand{\finprf}{\unskip\null\hfill$\square$\vskip 0.3cm}
\newcommand{\R}{\mathbb R}
\newcommand{\nrmu}[2]{\|#1\|_{L^{#2}(\mu)}}
\newcommand{\nrnu}[2]{\|#1\|_{L^{#2}(\nu)}}
\newcommand{\be}[1]{\begin{equation}\label{#1}}
\newcommand{\ee}{\end{equation}}
\newcommand{\eqn}[1]{(\ref {#1})}
\newcommand{\eps}{\varepsilon}
\newcommand{\vpi}{\pi}

\begin{document}
\begin{frontmatter}
\selectlanguage{english}
\title{Convex Sobolev inequalities and spectral gap}
\vspace{-2.6cm}
\selectlanguage{francais}
\title{In\'egalit\'es de Sobolev convexes et trou spectral}
\selectlanguage{english}
\author[authorlabel1]{Jean Dolbeault}
\ead{dolbeaul@ceremade.dauphine.fr}
\author[authorlabel1,authorlabel2]{Jean-Philippe Bartier}
\ead{bartier@math.uvsq.fr}
\address[authorlabel1]{Ceremade (UMR CNRS no. 7534), Universit\'e Paris-Dauphine, pl. de Lattre de Tassigny, 75775 Paris C\'edex~16, France}\address[authorlabel2]{Laboratoire de Math{\'e}mati\-ques Appliqu{\'e}es (UMR CNRS 7641), Universit\'e de Versailles, 45, av. des \'Etats-Unis, 78035 Versailles, France.}
\selectlanguage{english}
\begin{abstract} This note is devoted to the proof of convex Sobolev (or generalized Poincar\'e) inequalities which interpolate between spectral gap (or Poincar\'e) inequalities and logarithmic Sobolev inequalities. We extend to the whole family of convex Sobolev inequalities results which have recently been obtained by Cattiaux \cite{Cattiaux-2004} and Carlen and Loss \cite{MR2079069} for logarithmic Sobolev inequalities. Under local conditions on the density of the measure with respect to a reference measure, we prove that spectral gap inequalities imply all convex Sobolev inequalities with constants which are uniformly bounded in the limit approaching the logarithmic Sobolev inequalities. We recover the case of the logarithmic Sobolev inequalities as a special case.
\vskip 0.5\baselineskip
\selectlanguage{francais}
\noindent{\bf R\'esum\'e} Cette note est consacr\'ee \`a la preuve d'in\'egalit\'es de Sobolev convexes (ou in\'egalit\'es de Poincar\'e g\'en\'erali\-s\'ees) qui interpolent entre des in\'egalit\'es de trou spectral (ou de Poincar\'e) et des in\'egalit\'es de Sobolev logarithmiques. Nous \'etendons \`a la famille des in\'egalit\'es de Sobolev convexes toute enti\`ere des r\'esultats qui ont \'et\'e obtenus r\'ecemment par Cattiaux \cite{Cattiaux-2004} et Carlen et Loss \cite{MR2079069} pour des in\'egalit\'es de Sobolev logarithmiques. Sous des conditions locales sur la densit\'e de la mesure par rapport \`a une mesure de r\'ef\'erence, nous d\'emontrons que les in\'egalit\'es de trou spectral entra\^{\i}nent toutes les in\'egalit\'es de Sobolev convexes avec des constantes qui sont born\'ees uniform\'ement dans la limite qui approche les in\'egalit\'es de Sobolev logarithmiques. Nous retrouvons le cas des in\'egalit\'es de Sobolev logarithmiques comme un cas particulier.
\vskip 0.5\baselineskip
\end{abstract}
\end{frontmatter}
%
%
%

\section*{Version fran\c{c}aise abr\'eg\'ee}

Soit $\mu$ une mesure de probabilit\'e sur $\R^d$. On dira que $\mu$ admet une {\sl inegalit\'e de Sobolev logarithmique (tendue)\/} s'il existe une constante $C_1(\mu)$ telle que 
\be{Ineq:LSI1}
\iomu{u^2\,\log\Big(\frac{u^2}{\nrmu u2^2}\Big)}\le C_1(\mu)\,\nrmu{\nabla u}2^2\quad\forall\;u\in H^1(\mu)\;,
\ee
et une {\sl in\'egalit\'e de Poincar\'e,\/} ou encore de {\sl trou spectral,\/} s'il existe une constante $C_2(\mu) $ telle que
\be{Ineq:Poincare1}
\iomu{|u-\bar u|^2}\le C_2(\mu)\,\nrmu{\nabla u}2^2\quad\mbox{avec}\quad\bar u=\iomu u\;,\quad\forall\;u\in H^1(\mu)\;.
\ee
Si l'in\'egalit\'e \eqn{Ineq:LSI1} est v\'erifi\'ee, alors \eqn {Ineq:Poincare1} est aussi vraie avec $C_2(\mu)\le \frac 12\,C_1(\mu)$. La r\'eciproque est fausse en g\'en\'eral (consid\'erer $d\mu(x)=C\,\exp(-\vert x \vert^\alpha)$ avec $\alpha \in[1,2)$, voir \cite{MR1682772,Gentil-Guillin-Miclo} et \cite{MR1796718,MR2052235,BCR2004} pour plus de d\'etails). Cependant, Cattiaux dans \cite{Cattiaux-2004} puis Carlen et Loss dans \cite{MR2079069} ont donn\'e des conditions n\'ecessaires sur $\mu$ pour que \eqn{Ineq:LSI1} se d\'eduise de \eqn {Ineq:Poincare1}. Le but de cette note est d'am\'eliorer certains de ces r\'esultats en consid\'erant une famille d'in\'egalit\'es qui interpole entre \eqn{Ineq:LSI1} et~\eqn{Ineq:Poincare1}. A la suite de Beckner \cite{MR954373}, pour $p\in (1,2]$, on dira que $\mu$ v\'erifie une {\sl in\'egalit\'e de Poincar\'e g\'en\'eralis\'ee\/} s'il existe une constante positive finie $C_p(\mu)$ telle que
\be{Ineq:Beckner1}
\frac 1{p-1}\left[\,\iomu{|u|^2}-\Big(\iomu{|u|^{2/p}}\Big)^p\right]\le C_p(\mu)\iomu{|\nabla u|^2}\quad\forall\;u\in H^1(\mu)\;.
\ee
Le cas $p=2$ correspond \`a \eqn {Ineq:Poincare1} et dans la limite $p\rightarrow1$ , on retrouve \eqn{Ineq:LSI1} si $\liminf_{p\rightarrow 1}C_p(\mu)$ est finie. On peut montrer que
\be{Ineq:Encadrement-f}
\frac 2p\,C_2(\mu)\le C_p(\mu)\quad\forall\; p\in[1,2] \quad\mbox{et}\quad C_p(\mu)\le\frac 1{p-1}\, C_2(\mu)\quad\forall\; p\in(1,2]
\ee
(voir \cite{Arnold-Dolbeault,BDIK} pour plus de d\'etails). On en d\'eduit que si $C_1(\mu)$ est finie, alors $C_p(\mu)$ est aussi finie pour tout $p\in (1,2]$. Par contre, il n'est pas possible d'en d\'eduire une estimation de $C_1(\mu)$ sachant que $C_2(\mu)$ est finie. Dans la suite, $C_p(\mu)$ d\'esignera pour tout $p\in [1,2]$ la valeur optimale de la constante.

\medskip Gross a montr\'e dans \cite{MR0420249} que \eqn{Ineq:LSI1} est v\'erifi\'ee pour les mesures gaussiennes:
\[
\mu(x)=\nu_\sigma(x):=(2\pi\sigma^2)^{-d/2}\,e^{-\frac{|x|^2}{2\sigma^2}}\;,
\]
et, en utilisant les polyn\^omes d'Hermite, Beckner a \'etabli dans \cite{MR954373} que \eqn{Ineq:Beckner1} est aussi v\'erifi\'ee lorsque $\mu=\nu_\sigma$, pour tout $p\in (1,2)$, avec
\[
C_p(\nu_\sigma)=\frac 2p\,\sigma^2\,.
\]
La m\'ethode d'entropie -- production d'entropie de Bakry et Emery \cite{MR772092} permet de montrer \eqn{Ineq:Beckner1} dans le cas de mesures du type $\mu=e^{-V}$ lorsque $V$ est strictement convexe, voir \cite{MR1842428}. Pour cette raison, les in\'egalit\'es~\eqn{Ineq:Beckner} sont aussi appel\'ees {\sl in\'egalit\'es de Sobolev convexes.\/} On montre ainsi que
\[
C_p(e^{-V})\;\le\;\frac 2p\;\left[\inf_{\xi\in S^{d-1},\,x\in\R^d}\;(\,D^2V(x)\,\xi,\,\xi\,)\right]^{-1}\,.
\]
Cela conduit naturellement \`a rechercher des conditions suffisantes sur $V$ pour borner $C_1(\mu)$ en fonction de $C_2(\mu)$, ou, en d'autres termes, pour que l'in\'egalit\'e de Poincar\'e entra\^{\i}ne l'in\'egalit\'e de Sobolev logarithmique. Nous allons nous int\'eresser \`a des estimations de $C_p(\mu)$ pour tout $p\in (1,2)$, et, en prenant la limite $p\to 1$, retrouver et am\'eliorer les r\'esultats obtenus pour $p=1$ par Cattiaux dans \cite{Cattiaux-2004} puis par Carlen et Loss dans \cite{MR2079069}. Notre principal r\'esultat est un r\'esultat de perturbation pour les in\'egalit\'es de Sobolev convexes \eqn{Ineq:Beckner1} qui diff\`ere toutefois de la m\'ethode classique de Holley-Stroock~\cite{MR893137,MR1842428,Arnold-Dolbeault}.
\begin{theoreme}\label{Thm:Main1} Soit $p\in[1,2)$ et $p'=p/(p-1)$. Si $\mu$ et $\nu$ sont deux mesures de probabilit\'es de densit\'es respectives $e^{-V}$ et $e^{-W}$ par rapport \`a la mesure de Lebesgue telles que $C_p(\nu)$ et $C_2(\mu)$ soient finies et si $Z:=\frac 12(V-W)$ est une fonction de $L^{p'}(d\nu)$ telle que
\[
\inf_{x\in\R^d}\left(|\nabla Z|^2-\Delta Z+\nabla Z\cdot\nabla W\right)>-\infty\;,
\]
alors
\[
C_p(\mu)\le \mathcal C_p:=\frac 2p\,C_2(\mu)+\Big(\frac 2p-1\Big)\,\left[C_p(\nu)+C_2(\mu)\,\big(2\,\nrmu Z{p'}-m\,C_p(\nu)\big)_+\right]\,.
\]
\end{theoreme}

\medskip\noindent Par passage \`a la limite $p\to 1$, on obtient un r\'esultat pour les in\'egalit\'es de Sobolev logarithmiques~\eqn{Ineq:LSI1}.\smallskip
\begin{corollaire} Avec les m\^emes notations que ci-dessus, si les hypoth\`eses du Th\'eor\`eme \ref{Thm:Main1} sont v\'erifi\'ees uniform\'ement dans la limite $p\to 1$, alors $\mu$ v\'erifie l'in\'egalit\'e de Sobolev logarithmique \eqn{Ineq:LSI1} avec $C_1(\mu)\leq \liminf_{p\rightarrow 1} \mathcal C_p$.\end{corollaire} 

\medskip La preuve du Th\'eor\`eme \ref{Thm:Main1} consiste comme dans \cite{MR2079069} pour $p=1$ \`a \'etablir d'abord une in\'egalit\'e restreinte: \[
\frac{\iomu {|v|^2}-\Big(\iomu{|v|^{2/p}}\Big)^p}{(p-1)\iomu{|\nabla v|^2}}
\le\mathcal C^*_p\quad\forall\; v\in H^1(\mu)\quad\mbox{tel que}\quad\bar v=0\;.
\]
Ensuite on en d\'eduit le cas g\'en\'eral gr\^ace au
\begin{lemme}\label{Lem:Bakry-1} Soit $q\in[1,2]$. Pour toute fonction $u\in L^1\cap L^q(\mu)$, si $\bar u:=\iomu u$, alors
\[
\Big(\,\iomu{|u|^q}\,\Big)^{2/q}\ge |\bar u|^2+(q-1)\,\Big(\,\iomu{|u-\bar u|^q}\,\Big)^{2/q}\,.
\]\end{lemme}

\par\medskip\centerline{\rule{2cm}{0.2mm}}\medskip
\selectlanguage{english}
\setcounter{equation}{0}
\section{Introduction and main result}\label{Sec:Introduction}

Consider a probability measure $\mu$ on $\R^d$. We say that there is a {\sl (tight) logarithmic Sobolev inequality\/} associated to $\mu$ if there exists a finite constant $C_1(\mu)$ such that 
\be{Ineq:LSI}
\iomu{u^2\,\log\Big(\frac{u^2}{\nrmu u2^2}\Big)}\le C_1(\mu)\,\nrmu{\nabla u}2^2\quad\forall\;u\in H^1(\mu)\;,
\ee
and a {\sl Poincar\'e inequality\/} associated to $\mu$ if there exists a finite constant $C_2(\mu)$ such that 
\be{Ineq:Poincare}
\iomu{|u-\bar u|^2}\le C_2(\mu)\,\nrmu{\nabla u}2^2\quad\mbox{with}\quad\bar u=\iomu u\;,\quad\forall\;u\in H^1(\mu)\;.
\ee
This inequality is often called the {\sl spectral gap inequality\/} for the following reason. Consider in $H^1(\mu)$ the Rayleigh quotient $\nrmu{\nabla u}2^2/\nrmu u2^2$. The lowest critical value, zero, corresponds to constant functions, and the optimal value for $C_2(\mu)^{-1}$ is therefore associated with the second critical value: $ u-\bar u$ is the projection on the orthogonal of the constants with respect to the $L^2(\mu)$ norm. It is well known that if \eqn{Ineq:LSI} holds, then \eqn{Ineq:Poincare} is also true with 
\[
C_2(\mu)\le \frac 12\,C_1(\mu)\;.
\]
This is easily checked by writing $u=1+\eps\,v$, with $\bar v=0$, and by letting $\eps\to0$. The reverse implication is a much harder question, and not true in general. With no additional assumption, we may have $C_1(\mu)=+\infty$ and $C_2(\mu)<\infty$. An example of such a situation is given by $\mu(x)=\exp(-\vert x \vert ^\alpha)$ in $\R^d$ with $\alpha\in[1,2)$, see, {\sl e.g.,\/} \cite{MR1682772,Gentil-Guillin-Miclo} and \cite{MR1796718,MR2052235,BCR2004} for more details. Cattiaux in \cite{Cattiaux-2004}, and then Carlen and Loss in \cite{MR2079069} with more elementary tools, gave sufficient conditions on $\mu$ under which \eqn{Ineq:LSI} is a consequence of \eqn{Ineq:Poincare}. The goal of this note is to revisit some of these results by considering a family of inequalities which interpolate between \eqn{Ineq:LSI} and \eqn{Ineq:Poincare}. 

According to Beckner in \cite{MR954373}, we shall say that, for some $p\in (1,2]$, there is a {\sl generalized Poincar\'e inequality\/} associated to $\mu$ if there exists a positive constant $C_p(\mu)$ such that 
\be{Ineq:Beckner}
\frac 1{p-1}\left[\,\iomu{|u|^2}-\Big(\iomu{|u|^{2/p}}\Big)^p\right]\le C_p(\mu)\iomu{|\nabla u|^2}\quad\forall\;u\in H^1(\mu)\;.
\ee
Throughout this paper, we will assume that for any $p\in[1,2]$, $C_p(\mu)$ is the optimal constant. We will not consider ``defective'' logarithmic Sobolev inequality (see, {\sl e.g.,\/} \cite{Cattiaux-2004,MR1813804}) and will omit the word ``tight'' whenever we mention Inequality \eqn{Ineq:LSI}. The limit case $p=2$ corresponds to \eqn{Ineq:Poincare}, at least for nonnegative solutions. However, in the general case, \eqn{Ineq:Poincare} looks different of \eqn{Ineq:Beckner} in the limit case $p=2$. We indeed get
\[\label{Ineq:Beckner-p=2}
\iomu{|u|^2}-\Big(\iomu{|u|}\Big)^2\le C_2(\mu)\iomu{|\nabla u|^2}\quad\forall\;u\in H^1(\mu)
\]
in that case, which is equivalent to \eqn{Ineq:Poincare} only for nonnegative functions. However, if the inequality holds for a function $u-a=v\geq 0$, a straightforward computation shows that
\[
\iomu{|u-\bar u|^2}=\iomu{|(u-a)-(\bar u-a)|^2}=\iomu{|v-\bar v|^2}\le C_2(\mu)\iomu{|\nabla v|^2}=C_2(\mu)\iomu{|\nabla u|^2}\,,
\]
so that \eqn{Ineq:Poincare} holds for any $u\in H^1(\mu)$ such that $u_-\in L^\infty(\mu)$. By density, we extend it to any $u\in H^1(\mu)$: \eqn{Ineq:Beckner} with $p=2$ is therefore equivalent to \eqn{Ineq:Poincare}. On the other hand, by taking the limit $p\to 1$ in \eqn{Ineq:Beckner}, we find $C_1(\mu)\le\liminf_{p\to 1}C_p(\mu)$, which proves \eqn{Ineq:LSI} if the right hand side is finite. By considering again $u=1+\eps\,v$, with $\bar v=0$, in the limit $\eps\to0$, we get: 
\[
C_p(\mu)\ge \frac 2p\,C_2(\mu)\quad\forall\; p\in[1,2]\;.
\]
By H\"older's inequality, $(\iomu u)^2\le(\iomu{|u|^{2/p}})^p$ for any $p\in[1,2]$. As a consequence, for any $p\in(1,2]$,
\[
C_p(\mu)=\sup_{u\in H^1(\mu)}\kern -5pt\frac{\iomu {|u|^2}-\Big(\iomu{|u|^{2/p}}\Big)^p}{(p-1)\iomu{|\nabla u|^2}}\le \frac1{p-1}\,\sup_{u\in H^1(\mu)}\kern -5pt\frac{\iomu {|u|^2}-\Big(\iomu u\Big)^2}{\iomu{|\nabla u|^2}}=\frac{C_2(\mu)}{p-1}\;.
\]
We refer to \cite{Arnold-Dolbeault,BDIK} for more details. Summarizing, we know that 
\be{Ineq:Encadrement}
\frac 2p\,C_2(\mu)\le C_p(\mu)\quad\forall\; p\in[1,2] \quad\mbox{and}\quad C_p(\mu)\le\frac 1{p-1}\, C_2(\mu)\quad\forall\; p\in(1,2]\;.
\ee
It follows that if $C_1(\mu) <\infty$, then for all $p\in (1,2]$, $C_p(\mu)<\infty$. However, at this stage, it is clear that we have no estimate on $C_1(\mu)$ if we only know that $C_p(\mu)$ is finite for some $p\in (1,2]$. 

\medskip Inequality \eqn{Ineq:LSI} has been established by Gross in \cite{MR0420249} in the case of Gaussian measures:
\[
\mu(x)=\nu_\sigma(x):=(2\pi\sigma^2)^{-d/2}\,e^{-\frac{|x|^2}{2\sigma^2}}\;,
\]
and using Hermite polynomials Beckner in \cite{MR954373} proved that \eqn{Ineq:Beckner} holds with
\[
C_p(\nu_\sigma)=\frac 2p\,\sigma^2\,.
\]
An alternative method based on the entropy -- entropy production method of Bakry and Emery \cite{MR772092} has been adapted in \cite{MR1842428} to prove \eqn{Ineq:Beckner} in more general situations which for instance cover the case of measures $\mu=e^{-V}$ for some strictly convex function $V$. For this reason, the family of inequalities \eqn{Ineq:Beckner} has been called {\sl convex Sobolev inequalities.\/} The entropy -- entropy production method gives an upper bound on the best constant:
\[
C_p(e^{-V})\;\le\;\frac 2p\;\left[\inf_{\xi\in S^{d-1},\,x\in\R^d}\;(\,D^2V(x)\,\xi,\,\xi\,)\right]^{-1}\;=:\;\frac 2{p\,\lambda_1}\;.
\]
This shows that at least in some circumstances, the bounds in \eqn{Ineq:Encadrement} are not optimal. As already mentioned, Cattiaux in \cite{Cattiaux-2004} and then Carlen and Loss in \cite{MR2079069} gave sufficient conditions on $V$ to bound $C_1(\mu)$ in terms of $C_2(\mu)$, or, in other words, to deduce logarithmic Sobolev inequalities from spectral gap inequalities. Our purpose is to extend these results to $C_p(\mu)$ for any $p\in (1,2)$, and recover and improve their results by deriving uniform estimates in the limit $p\to 1$. 

Note for completeness that improvements of several types have been obtained, for instance by considering $L^\infty$ perturbations of $V$ based on Holley-Stroock type estimates \cite{MR893137,MR1842428,Arnold-Dolbeault}. This allows to relax the strict convexity condition on $V$. One can also refine the entropy -- entropy production method \cite{Arnold-Dolbeault}, thus giving for instance the improved inequality
\[
\Big(\frac p{p-1}\Big)^2\left[\,\iomu{|u|^2}-\Big(\iomu{|u|^{2/p}}\Big)^{2(p-1)}\Big(\iomu{|u|^2}\Big)^{\frac 2p-1}\right]\le \frac 4{\lambda_1}\iomu{|\nabla u|^2}\quad\forall\;u\in H^1(\mu)\;.
\]
Although it differs in nature from the Holley and Stroock perturbation lemma, our main result can be seen as a perturbation result as well. It applies to convex Sobolev inequalities \eqn{Ineq:Beckner}.\smallskip
\begin{theorem}\label{Thm:Main} Let $p\in [1,2)$. Let $\mu$ and $\nu$ be two probability measures with respective densities $e^{-V}$ and $e^{-W}$ relatively to Lebesgue's measure such that, for some $p \in(1,2]$, $C_p(\nu)$ and $C_2(\mu)$ are finite. Assume that 
\[
Z:=\frac 12(V-W)\in L^{p'}(d\nu)\quad\mbox{and}\quad m:=\inf_{x\in\R^d}\delta(x)>-\infty\;,
\]
where $\delta:=|\nabla Z|^2-\Delta Z+\nabla Z\cdot\nabla W$, and define $\mathcal C^*_p:=C_p(\nu)+C_2(\mu)\,\big(2\,\nrmu Z{p'}-m\,C_p(\nu)\big)_+$. Then we have
\[
C_p(\mu)\le \mathcal C_p:=\frac 2p\,C_2(\mu)+\Big(\frac 2p-1\Big)\,\mathcal C^*_p\;.
\]
\end{theorem}

\medskip\noindent  We denote by $p'=p/(p-1)\in[2,\infty]$ the H\"older conjugate of $p\in[1,2]$. By relative density, we simply mean, {\sl e.g.,\/} $d\mu(x)=e^{-V(x)}\,dx$. With these notations, $\mu=e^{-2Z}\,\nu$. Taking the limit $p\to 1$, we obtain a result analoguous to Theorem \ref{Thm:Main} for the logarithmic Sobolev inequalities \eqn{Ineq:LSI}.\smallskip
\begin{corollary} With the above notations, if the assumptions of Theorem \ref{Thm:Main} hold uniformly in the limit $p\to 1$ and if $\liminf_{p\rightarrow 1} \mathcal C_p$ is finite, then the logarithmic Sobolev inequality \eqn{Ineq:LSI} associated to $\mu$ holds with $C_1(\mu)\leq \liminf_{p\rightarrow 1} \mathcal C_p$. 
\end{corollary} 

\section{Proof of the main result}\label{Sec:Proof}

As in \cite{MR2079069}, we first prove Theorem \ref{Thm:Main} in the {\sl restricted case\/} which corresponds to $\bar u=0$ and then extend it to the {\sl unrestricted case.\/} 
\begin{lemma}\label{Lem:restricted} Under the assumptions of Theorem~\ref{Thm:Main}, 
\[
\sup_{v\in H^1(\mu),\,\bar v=0}\frac{\iomu {|v|^2}-\Big(\iomu{|v|^{2/p}}\Big)^p}{(p-1)\iomu{|\nabla v|^2}}
\le\mathcal C^*_p
\;.
\]\end{lemma}
\proof Define
\[
\mathcal A(t):=\nrmu{\nabla v}2^2-\frac{t}{(p-1)\,C_p(\nu)}\left[\,\iomu{|v|^2}-\Big(\iomu{|v|^{2/p}}\Big)^p\right]\,.
\]
Proving that for some $t>0$, $\mathcal A (t)\ge 0$ for any $v$ in $H^1(\mu)$ with $\bar v=0$ is equivalent to the result of Theorem~\ref{Thm:Main} in the case $\bar u=0$, {\sl i.e.\/} $u=v$, the so-called {\sl restricted case\/} in \cite{MR2079069}. Let us write
\[
\mathcal A(t)=\mathrm{(I)}+\mathrm{(II)}+\mathrm{(III)}
\]
with
\[\begin{array}{l}
\mathrm{(I)}=(1-t)\iomu{|\nabla v|^2}\;,\\
\mathrm{(II)}=t\iomu{|\nabla v|^2}\;,\\
\mathrm{(III)}=\frac{-t}{(p-1)\,C_p(\nu)}\left[\,\iomu{|v|^2}-\Big(\iomu{|v|^{2/p}}\Big)^p\right]\;.\\
\end{array}\]
Let $v=g\,e^Z$:
\[
\iomu{|v|^2}=\iomuz{|g|^2}\quad\mbox{and}\quad\iomu{|\nabla v|^2}=\iomuz{|\nabla g|^2}+\iomuz{\delta\,|g|^2}\;.
\]
Using the spectral gap assumption on $\mu$, we get
\[
\mathrm{(I)}\ge \frac{1-t}{C_2(\mu)}\iomu{|v|^2} = \frac{1-t}{C_2(\mu)}\iomuz{|g|^2}\;.
\]
Using the fact that \eqn{Ineq:Beckner} holds for $\nu$ and the above expression of $\iomu{|\nabla v|^2}$, we obtain
\[
\mathrm{(II)} \ge \frac t{(p-1)\,C_p(\nu)}\left(\,\iomuz{|g|^2}-\Big(\iomuz{|g|^{2/p}}\Big)^p\right)+t\iomuz{\delta\,|g|^2}\,.
\]
As for the last term, we can write it as
\[
\mathrm{(III)}=\frac t{(p-1)\,C_p(\nu)}\left(\Big(\iomu{|v|^{2/p}}\Big)^p-\iomuz{|g|^2}\right)\,.
\]
Collecting these estimates, we have
\[
\mathcal A(t) \ge \iomu{\Big(\frac{(1-t)}{C_2(\mu)}+t\,\delta\Big)|g|^2}+\frac{\mathcal B\,t}{(p-1)\,C_p(\nu)} \;,\quad \mathcal B:=\Big(\iomu{|v|^{2/p}}\Big)^p-\Big(\iomuz{|g|^{2/p}}\Big)^p.
\]
Let $\d\vpi:=\frac{|g|^{2/p}}{\iomuz{|g|^{2/p}}}d\nu$. By Jensen's inequality applied to the convex function $t\mapsto e^{-t}$, we get
\[
\frac{\iomu{|v|^{2/p}}}{\iomuz{|g|^{2/p}}}=\frac{\iomuz{|g|^{2/p}e^{-2(1-\frac 1p)Z}}}{\iomuz{|g|^{2/p}}}=\int e^{-2(1-\frac 1p)Z}\,d\vpi\ge \exp\Big[-2\Big(1-\frac 1p\Big)\int Z\,d\vpi\Big]\,.
\]
Using the lower estimate $e^{-t}\ge 1-t$, we infer that 
\[
\left(\frac{\iomu{|v|^{2/p}}}{\iomuz{|g|^{2/p}}}\right)^p=e^{-2(p-1)\int Z\,d\vpi}\ge 1-2(p-1)\int Z\,d\vpi\;,
\]
\[
\Big(\iomu{|v|^{2/p}}\Big)^p-\Big(\iomuz{|g|^{2/p}}\Big)^p\ge-2(p-1)\iomuz{Z\,|g|^{2/p}}\Big(\iomuz{|g|^{2/p}}\Big)^{p-1}\;.
\]
By H\"older's inequality, we have
\[
\Big(\iomuz{|g|^{2/p}}\Big)^{p-1}\le \Big(\iomuz{|g|^2}\Big)^{1-1/p}\quad\mbox{and}\quad\iomuz{Z\,|g|^{2/p}}\le\Big(\iomuz{|g|^2}\Big)^{1/p}\;\nrnu Z{p'}\,,
\]
\[
\mathcal B\ge -2\,(p-1)\,\nrmu Z{p'}\iomuz{|g|^{2}}\;.
\]
Altogether, we get
\[
\mathcal A(t)\ge\left[\frac{1-t}{C_2(\mu)}+t\left(m-\frac{2\,\nrmu Z{p'}}{C_p(\nu)}\right)\right]\iomuz{ |g|^{2}}\;.
\]
This proves that $\mathcal A(t)\ge 0$ for any $t\in (0,t^*]$ with
\[
t^*\!:=1\;\mbox{if}\; m-\frac 1{C_2(\mu)}-\frac{2\,\nrmu Z{p'}}{C_p(\nu)}\ge 0\;,\quad t^*\!:=\left[1+\frac{C_2(\mu)}{C_p(\nu)}\,\Big(2\,\nrmu Z{p'}-m\,C_p(\nu)\Big)\right]^{-1}\mbox{otherwise}\,.
\]
This ends the proof with $\mathcal C_p=C_p(\nu)/t^*$. \finprf

The general case $\bar u\neq 0$ in Theorem \ref{Thm:Main} is a consequence of the following estimate, which is the counterpart for $p<2$ of a Lemma given in \cite{Bakry92} for $p>2$ (see Remark 2 below).
\begin{lemma}\label{Lem:Bakry} Let $q\in[1,2]$. For any function $u\in L^1\cap L^q(\mu)$, if $\bar u:=\iomu u$, then
\[
\Big(\,\iomu{|u|^q}\,\Big)^{2/q}\ge |\bar u|^2+(q-1)\,\Big(\,\iomu{|u-\bar u|^q}\,\Big)^{2/q}\,.
\]\end{lemma}
\proof Let $v:=u-\bar u$, $\phi(t):=\left(\,\iomu{|\bar u+t\,v|^q}\,\right)^{2/q}$. We may notice that $\phi(0)=|\bar u|^2$, $\phi'(0)=0$, $\phi(1)=(\,\iomu{|u|^q}\,)^{2/q}$ and 
\[
\frac 12\,\phi''(t)=(2-q)\Big(\iomu{|w|^q}\Big)^{\frac 2q-2}\Big(\iomu{|w|^{q-2}\,w\,v}\,\Big)^2+(q-1)\Big(\iomu{|w|^q}\,\Big)^{\frac 2q-1}\iomu{|w|^{q-2}v^2}
\]
with $w:=\bar u+t\,v$. The first term of the right hand side is nonnegative. As for the second one, we may use H\"older's inequality:
\[
\Big(\,\iomu{|v|^q}\,\Big)^{\frac 2q}=\Big(\,\iomu{|w|^{\frac q2(2-q)}\cdot|v|^q\,|w|^{\frac q2(q-2)}}\,\Big)^{\frac 2q}\le \Big(\,\iomu{|w|^q}\,\Big)^{\frac 2q-1}\cdot\iomu{|w|^{q-2}\,|v|^2}\;.
\]
Thus we get: $\frac 12\,\phi''(t)\ge (q-1)\,(\iomu{|v|^q})^{2/q}$, which proves that $\phi(1)\ge\phi(0)+(q-1)\,(\,\iomu{|v|^q}\,)^{2/q}$ and completes the proof. \finprf

\noindent{\sl Proof of Theorem \ref{Thm:Main}.\/} Let $v:=u-\bar u$ and apply Lemma~\ref{Lem:Bakry} with $q=\frac 2p\in[1,2)$. Since $\iomu{|u|^2}-|\bar u|^2=\iomu{|u-\bar u|^2}=\iomu{|v|^2}$, we can write
\begin{eqnarray*}
\iomu{|u|^2}-\Big(\iomu{|u|^{2/p}}\Big)^p&\le&\iomu{|u|^2}-|\bar u|^2-\Big(\frac 2p-1\Big)\left(\iomu{|u-\bar u|^{\frac 2p}}\right)^p\\
&&=\iomu{|v|^2}-\Big(\frac 2p-1\Big)\left(\iomu{|v|^{\frac 2p}}\right)^p\\
&&=2\,\frac{p-1}p\,\iomu{|v|^2}+\frac{2-p}p\,\left[\,\iomu{|v|^2}-\Big(\iomu{|v|^{2/p}}\Big)^p\,\right]\,.
\end{eqnarray*}
We can then apply \eqn{Ineq:Poincare} and Lemma \ref{Lem:restricted}, and the result holds with $\mathcal C_p =\frac 2p\,C_2(\mu)+(\frac 2p-1)\,\mathcal C^*_p$. \finprf
\begin{remark} -- To deduce the unrestricted inequality from the restricted inequality, Carlen and Loss in \cite{MR2079069} use the following inequality:
\[
\iomu{|u|^q}\le |\bar u|^q+\frac 12\,q\,(q-1)\,\nrmu uq^{q-2}\,\nrmu vq^2\quad\forall u\in L^q(\mu)\,,\quad v=u-\bar u\;\quad\forall\;q\in [2,\infty)\;.
\]
The proof is essentially the same as for Lemma~\ref{Lem:Bakry}. We can also write a similar result for $q\leq 2$:
\[
\iomu{|u|^q}\ge |\bar u|^q+\frac 12\,q\,(q-1)\,\nrmu uq^{q-2}\,\nrmu vq^2\quad\forall u\in L^q(\mu)\,,\quad v=u-\bar u\;\quad\forall\;q\in (1,2]\;.
\]\end{remark}
\begin{remark} -- In the case $q>2$, according to \cite{Bakry92}, the following result holds:
\[
\Big(\,\iomu{|u|^q}\,\Big)^{2/q}\le |\bar u|^2+(q-1)\,\Big(\,\iomu{|u-\bar u|^q}\,\Big)^{2/q}\,.
\]
\end{remark}
\begin{remark} -- For evident reasons, $\mathcal C_p^*\leq \mathcal C_p$: to the restricted case corresponds an improved inequality, stated in Lemma~\ref{Lem:restricted}. On the other hand, from \eqn{Ineq:Encadrement} and Theorem  \ref{Thm:Main}, we obtain
\[
0\le \Big(\frac 2p-1\Big)\,C_2(\mu)\le C_p(\mu)-C_2(\mu)\le\mathcal C_p-C_2(\mu)\le \frac{2-p}p\,\left(C_2(\mu)-\mathcal C_p^*\right)\;.
\]
This means that for any $p\in (1,2)$, under the assumptions of Theorem \ref{Thm:Main},
\[
\mathcal C_p^*\le C_2(\mu)\le\mathcal C_p\;.
\]
\end{remark}

\section{Application to the euclidean space}\label{Sec:Application}

To compare our results with those of \cite{MR2079069}, we can state a result for generalized Poincar\'e inequalities corresponding to Gaussian weights, {\sl i.e.\/} $\nu=\nu_\sigma$, and recover in the limit $p\to 1$ the logarithmic Sobolev inequality. We can optimize the choice of $W(x)=|x|^2/(2\sigma^2)$ and cover, for instance, all harmonic potentials, which was not the case in \cite{MR2079069}. This freedom in the choice of the parameter $\sigma$ corresponds to the scaling invariance in the Euclidean space, which is however not so easy to write in the case of generalized Poincar\'e inequalities. In the case where $\nu=\nu_\sigma$, it is known that inequality (\ref{Ineq:Beckner}) holds. Theorem \ref{Thm:Main} becomes\smallskip
\begin{corollary} Let $\nu=e^{-V}$ a probabality measure. If there exists $\sigma \in(0,\infty)$ such that
\[ 
V-\frac{\vert x \vert^2}{2 \sigma^2}\in L^{p'}(d\nu_\sigma)\quad \mbox{and}\quad\inf_{x\in\R^d} \left(\left\vert\nabla V \right\vert ^2-2\,\Delta V -\frac {|x|^2}{\sigma^4}\right)> -\infty \; ,
\]
then Inequality \eqn{Ineq:LSI} holds.
\end{corollary}

\bibliographystyle{siam}
\bibliography{References}

\begin{thebibliography}{10}

\bibitem{Arnold-Dolbeault}
{\sc A.~Arnold and J.~Dolbeault}, {\em Refined convex {S}obolev inequalities},
  tech. rep., Ceremade, preprint no. 0431, 2004.

\bibitem{MR1842428}
{\sc A.~Arnold, P.~Markowich, G.~Toscani, and A.~Unterreiter}, {\em On convex
  {S}obolev inequalities and the rate of convergence to equilibrium for
  {F}okker-{P}lanck type equations}, Comm. Partial Differential Equations, 26
  (2001), pp.~43--100.

\bibitem{Bakry92}
{\sc D.~Bakry}, {\em L'hypercontractivit\'e et son utilisation en th\'eorie des
  semigroupes}, in Lectures on probability theory (Saint-Flour, 1992),
  vol.~1581 of Lecture Notes in Math., Springer, Berlin, 1994, pp.~1--114.

\bibitem{MR772092}
{\sc D.~Bakry and M.~{\'E}mery}, {\em Hypercontractivit\'e de semi-groupes de
  diffusion}, C. R. Acad. Sci. Paris S\'er. I Math., 299 (1984), pp.~775--778.

\bibitem{BCR2004}
{\sc F.~Barthe, P.~Cattiaux, and C.~Roberto}, {\em Interpolated inequalities
  between exponential and {G}aussian, {O}rlicz hypercontractivity and
  isoperimetry}, tech. rep., HAL Preprint no. 00002203, 2004.

\bibitem{MR2052235}
{\sc F.~Barthe and C.~Roberto}, {\em Sobolev inequalities for probability
  measures on the real line}, Studia Math., 159 (2003), pp.~481--497.
\newblock Dedicated to Professor Aleksander Pe\l czy\'nski on the occasion of
  his 70th birthday (Polish).

\bibitem{BDIK}
{\sc J.-P. Bartier, J.~Dolbeault, R.~Illner, and M.~Kowalczyk}, {\em A
  qualitative study of linear drift-diffusion equations with time-dependent or
  vanishing coefficients}.
\newblock 2005.

\bibitem{MR954373}
{\sc W.~Beckner}, {\em A generalized {P}oincar\'e inequality for {G}aussian
  measures}, Proc. Amer. Math. Soc., 105 (1989), pp.~397--400.

\bibitem{MR1682772}
{\sc S.~G. Bobkov and F.~G{\"o}tze}, {\em Exponential integrability and
  transportation cost related to logarithmic {S}obolev inequalities}, J. Funct.
  Anal., 163 (1999), pp.~1--28.

\bibitem{MR2079069}
{\sc E.~Carlen and M.~Loss}, {\em Logarithmic {S}obolev inequalities and
  spectral gaps}, in Recent advances in the theory and applications of mass
  transport, vol.~353 of Contemp. Math., Amer. Math. Soc., Providence, RI,
  2004, pp.~53--60.

\bibitem{Cattiaux-2004}
{\sc P.~Cattiaux}, {\em Hypercontractivity for perturbed diffusion
  semi-groups}, Preprint,  (2004).

\bibitem{Gentil-Guillin-Miclo}
{\sc I.~Gentil, A.~Guillin, and L.~Miclo}, {\em Modified logarithmic {S}obolev
  inequalities and transportation inequalities}, tech. rep., Ceremade, preprint
  no. 04-29, 2004.

\bibitem{MR0420249}
{\sc L.~Gross}, {\em Logarithmic {S}obolev inequalities}, Amer. J. Math., 97
  (1975), pp.~1061--1083.

\bibitem{MR893137}
{\sc R.~Holley and D.~Stroock}, {\em Logarithmic {S}obolev inequalities and
  stochastic {I}sing models}, J. Statist. Phys., 46 (1987), pp.~1159--1194.

\bibitem{MR1796718}
{\sc R.~Lata{\l}a and K.~Oleszkiewicz}, {\em Between {S}obolev and
  {P}oincar\'e}, in Geometric aspects of functional analysis, vol.~1745 of
  Lecture Notes in Math., Springer, Berlin, 2000, pp.~147--168.

\bibitem{MR1813804}
{\sc M.~Ledoux}, {\em The geometry of {M}arkov diffusion generators}, Ann. Fac.
  Sci. Toulouse Math. (6), 9 (2000), pp.~305--366.
\newblock Probability theory.

\end{thebibliography}
\smallskip\noindent{\sl \copyright~2005 by the authors. This paper may be reproduced, in its entirety, for non-commercial purposes.}
\end{document}